\documentclass[12pt]{amsart}
\usepackage[centertags]{amsmath}

\usepackage{amsfonts}

\usepackage{amssymb}

\title{Minimal elementary end extensions}
\author{James H. Schmerl}

\date{\today}
\def\into{\longrightarrow}

\def\harp{\hspace{-4pt} \upharpoonright \hspace{-4pt}}

\def\pa{{\sf PA}}

\def\wkl{{\sf WKL}_0}
\def\aca{{\sf ACA}_0}

\def\wklsd{\wkl^*}

\newcommand{\KK}{{\mathcal K}}

\newcommand{\MM}{{\mathcal M}}

\newcommand{\NN}{{\mathcal N}}

\def\lhdeq{\trianglelefteq}

\def\rhd{\vartriangleright}

 \DeclareMathOperator{\cf}{cf}
 \DeclareMathOperator{\cod}{Cod}

  \DeclareMathOperator{\Def}{Def}
   \DeclareMathOperator{\den}{den}

  \DeclareMathOperator{\Lt}{Lt}

\begin{document}
\maketitle

\begin{abstract} Suppose that $\MM \models \pa$ and ${\mathfrak X} \subseteq {\mathcal P}(M)$. 
If $\MM$ has a {\it finitely generated} elementary end extension $\NN \succ_{\sf end} \MM$ such that 
$\{X \cap M : X \in \Def(\NN)\} = {\mathfrak X}$, then there is such an $\NN$ that is, in addition, a {\it minimal} extension of $\MM$ 
iff every subset of $M$ that is $\Pi_1^0$-definable in $(\MM, {\mathfrak X})$ is the countable union of $\Sigma_1^0$-definable sets. 
  \end{abstract}

\smallskip


The fundamental theorem of MacDowell \& Specker \cite{ms} states that every model of Peano Arithmetic ({\sf PA}) has an elementary end extension. Since it first appeared more than half a century ago, this theorem has been frequently refined, generalized, modified and applied. 
This note  is in that spirit.

 For a model $\MM \models \pa$, we let $\Def(\MM)$ be the set of all parametrically definable subsets of $M$. If $\MM \prec \NN$, then the set of subsets of $M$ coded in $\NN$ is $\cod(\NN / \MM) = \{X \cap M : X \in \Def(N)\}$. An elementary extension $\NN \succ \MM$ is conservative iff $ \cod(\NN / \MM) = \Def(\MM)$. All conservative extensions are end extensions. Phillips \cite{phil},
 introducing the notion of a conservative extension,  improved the MacDowell-Specker Theorem by observing that every model of \pa\ has a conservative  extension. This was recently generalized by the following theorem. 

\bigskip

{\sc Theorem 1}: (\cite[Theorem~4]{code}) {\em  If $\MM \models \pa$  and  
${\mathfrak X} \subseteq {\mathcal P}(M)$, then the following are equivalent$:$

\begin{itemize}

\item[$(a)$] ${\mathfrak X}$ is countably generated, $\Def(\MM) \subseteq {\mathfrak X}$ and  
$(\MM, {\mathfrak X}) \models  \wklsd$. 

\item[$(b)$]
There is a countably generated extension   $\NN \succ_{\sf end} \MM$ such that 
${\mathfrak X} = \cod(\NN / \MM)$.

\item[$(c)$] There is a finitely generated extension  $\NN \succ_{\sf end} \MM$ such that \\
${\mathfrak X} = \cod(\NN / \MM)$.

\end{itemize}}

\bigskip

The theory $\wklsd$ of second-order arithmetic was introduced by Simpson \& Smith  \cite{ss} and is, roughly speaking, $\wkl$ with $\Sigma^0_1$-{\sf IND} replaced by $\Sigma_0^0$-{\sf IND}. We say that ${\mathfrak X}$ is countably generated if there is a countable ${\mathfrak X}_0 \subseteq {\mathfrak X}$ such that every set in ${\mathfrak X}$ is $\Delta_1^0$-definable in $(\MM, {\mathfrak X}_0)$.   It should be noted that  an immediate consequence  
 of \cite[Theorem~4.8]{ss} is that $(b)$ implies that $(\MM, {\mathfrak X}) \models \wklsd$. 

Improving the MacDowell-Specker Theorem  in another direction, Gaifman \cite{g} proved that every model of  
 \pa\ has a minimal elementary end extension. Phillips \cite{phil}  observed that this result could also be improved:    every model of \pa\ has a conservative 
 minimal  extension. He also proved  \cite{phil2} that, on the other hand,   the standard model  has a nonconservative minimal elementary extension.  These results were improved for countable models 
 by the following theorem. 

\bigskip

{\sc Theorem 2}: (\cite[Theorem~5]{code}) {\em  If $\MM \models \pa$ and  
${\mathfrak X} \subseteq {\mathcal P}(M)$, then the following are equivalent$:$

\begin{itemize}

\item[$(a)$] $\MM$ and ${\mathfrak X}$ are countable, $\Def(\MM) \subseteq {\mathfrak X}$ and  
$(\MM, {\mathfrak X}) \models  \wklsd$.

\item[$(b)$]
There is a countable  $\NN \succ_{\sf end} \MM$ such that 
${\mathfrak X} = \cod(\NN / \MM)$.

\item[$(c)$] There is a superminimal extension  $\NN \succ_{\sf end} \MM$ such that  \\
${\mathfrak X} = \cod(\NN / \MM)$.

\end{itemize}}

\bigskip

Recall that the extension $\NN \succ \MM$ is minimal iff there is no ${\KK}$ such that $\MM \prec \KK \prec \NN$, and it is superminimal iff 
 there is no ${\KK}$ such that $\MM \not\succcurlyeq \KK \prec \NN$. 
 Knight \cite{kn} proved that every countable model of \pa\ has a superminimal extension, and then this result was  extended  to conservative superminimal extensions (\cite[Cor.~2.2.12]{ksbook}).

Every superminimal extension is minimal, and every minimal extension is finitely generated. If we weaken the implications $(c) \Longrightarrow (a)$ of Theorem~1 by replacing ``finitely generated'' in (c) by ``minimal'' and   also weaken the implication $(a) \Longrightarrow (c)$ of Theorem~2 by replacing ``superminimal'' in (c) by ``minimal'', then  there is a wide gap between these two weakened implications. Question~6 of \cite{code} asked about the possibility of narrowing this gap. Theorem~3, the principal result of this paper,  answers this question by completely characterizing such $(\MM, {\mathfrak X})$.

\bigskip

{\sc Theorem 3}: {\em  If $\MM \models \pa$ and  
${\mathfrak X} \subseteq {\mathcal P}(M)$, then the following are equivalent$:$

\begin{itemize}

\item[$(a)$]  There is a countably generated extension  $\NN \succ_{\sf end} \MM$ such that  
${\mathfrak X} = \cod(\NN / \MM)$,  and every set that is $\Pi^0_1$-definable  in $(\MM, {\mathfrak X})$ is the union of   countably many $\Sigma_1^0$-definable sets.

\item[$(b)$] There is a minimal extension  $\NN \succ_{\sf end} \MM$ such that  \\
${\mathfrak X} = \cod(\NN / \MM)$.

\end{itemize}}

\bigskip

The additional condition  in $(a)$ holds whenever ${\mathfrak X}$ is countable. Also, it is implied by $\aca$, which is equivalent to: every 
$\Pi_1^0$-definable  set is $\Sigma_1^0$-definable. Thus,  the following corollary ensues.

\bigskip

{\sc Corollary 4}: {\em If   $({\mathcal M}, {\mathfrak X}) \models \aca$ and 
there is a countably generated extension $\NN \succ_{\sf end} \MM$ such that $\cod(\NN / \MM) = {\mathfrak X}$, then there is a minimal extension $\NN \succ_{\sf end} \MM$ such that $\cod(\NN / \MM) = {\mathfrak X}$.}

\bigskip

 $\Sigma_1^0$-{\sf IND} is equivalent to: every bounded $\Pi_1^0$-definable set is $\Sigma_1^0$-defin-able. If $(\MM, {\mathfrak X}) \models \Sigma_1^0$-${\sf IND}$ and $\MM$ has countable cofinality, then every $\Pi_1^0$-definable set is the union of countably many $\Sigma_1^0$-definable sets. Thus,   Corollary~4 can be improved when restricted 
to models of countable cofinality.

\bigskip

{\sc Corollary 5}: {\em If   $({\mathcal M}, {\mathfrak X}) \models \wkl$, $\cf(\MM) = \aleph_0$  and 
there is a countably generated extension $\NN \succ_{\sf end} \MM$ such that $\cod(\NN / \MM) = {\mathfrak X}$, then there is a minimal extension $\NN \succ_{\sf end} \MM$ such that \mbox{$\cod(\NN / \MM) = {\mathfrak X}$}.}

\bigskip

There are five numbered sections following this introduction. 
The purpose of  \S1 is to repeat some of the definitions from \cite{code} that are needed here. In particular, all definitions needed to understand Theorems~1--3 are given. 
The proof of a theorem that is stronger than the $(b) \Longrightarrow (a)$ half of Theorem~3 is proved in \S2. 
Since the proof of the other half of 
Theorem~3 depends on the proof of Theorem~1 as presented in \cite{code}, an outline of that proof is given in \S3.  
 Theorem~3 is proved in \S4. Finally, in \S5, we show with examples that the additional condition 
 in (a) of Theorem~3 is needed.

\bigskip


{\bf \S1.\@ Definitions.}  Let ${\mathcal L}_{\pa} = \{+,\times,0,1, \leq\}$ be the usual language appropriate for \pa. 
  Models of \pa\ will be denoted exclusively by script letters such as $\MM, \NN$, etc., and their respective  universes by the corresponding roman letters $M, N$, etc. A subset $R \subseteq M^n$ is $\MM$-definable if it is definable in $\MM$ possibly with parameters. 
  Thus, $\Def(\MM)$ is the set of all $\MM$-definable  subsets of $M$. 
  If $\MM \models \pa$ and $X \subseteq M$, then ${\mathcal L}_\pa(X) = {\mathcal L}_\pa \cup X$.

  Suppose that $\MM \prec \NN$. The definitions of when $\NN$ is a conservative, minimal or 
  superminimal extension of $\MM$ is given in the introduction. A subset $A \subseteq N$ {\it generates} 
$\NN$ {\it over} $\MM$ if whenever $\MM \preccurlyeq \KK \preccurlyeq \NN$ and 
$A \subseteq K$, then $\KK = \NN$. We say that  $\NN$ is a {\it countably} (respectively, {\it finitely}) {\it generated} extension of $\MM$ if there 
is a countable (respectively, finite) set $A \subseteq N$  that   generates $\NN$ over $\MM$. 

Suppose $\MM \models \pa$. Each $s \in M$ codes the $\MM$-finite sequence 
$\langle (s)_i : i < \ell(s) \rangle$, where $\ell(s)$, the {\it length} of $s$, is  the largest $y$ 
such that $2^y \leq s+1$ and whenever $i < \ell(s)$, then $(s)_i$ is the $i$-th digit in the binary expansion $s+1$. 
Define $\lhd$ on $M$ so that if $s,t \in M$, then $s \lhd t$ iff the sequence coded by $s$ is a proper initial segment of the one coded by $t$. In this way, we get the full binary $\MM$-tree 
$(M, \lhdeq)$, which we suggestively denote by~$2^{<M}$. An $\MM$-{\it tree} $T$ is a subset of $2^{<M}$ such that whenever $s \lhd t \in T$, then $s \in T$. 
A subset $B \subseteq T$ is a {\it branch} (of $T$) if it is an unbounded maximal linearly ordered subset. 
We will say that the branch $B$ {\it indicates} $X$ if 
$$
X = \{k \in M :  {\mbox{ for some }} s \in B,\ \MM \models \ell(s) > k \wedge (s)_k = 0  \}
$$
or, equivalently, 
$$
X = \{k \in M :  {\mbox{ for all }} s \in B,\ \MM \models \ell(s) > k \wedge (s)_k = 0 \}.
$$

 All models of second-order arithmetic encountered here 
 will have the form $(\MM, {\mathfrak X})$, 
 where $\MM \models \pa$ and ${\mathfrak X} \subseteq {\mathcal P}(M)$. 
 Consult Simpson's book \cite{simp} as a reference. Some schemes of second-order arithmetic will be needed. 
  We let 
$\wklsd = \Delta_1^0$-${\sf CA}_0 +  \Sigma_0^0$-{\sf IND} $+ {\sf WKL}$,  
 where ${\sf WKL}$ is Weak K\"{o}nig's Lemma, which asserts that every unbounded $\MM$-tree in ${ \mathfrak X}$ has a branch in ${\mathfrak X}$.  As usual, $\wkl = \wklsd + \Sigma_1^0$-{\sf IND} 
 and $\aca = \wkl + \Sigma^0_1$-${\sf CA}$.

  If ${\mathfrak X}_0 \subseteq 
{\mathfrak X} \subseteq   {\mathcal P}(M)$, 
then we say that ${\mathfrak X}_0$ {\it generates} ${\mathfrak X}$ if whenever 
${\mathfrak X}_0 \subseteq {\mathfrak X}_1 \subseteq {\mathfrak X}$, 
 and 
$(\MM, {\mathfrak X}_1) \models \Delta^0_1$-${\sf CA}$, then ${\mathfrak X}_1 = {\mathfrak X}$.
We say that  ${\mathfrak X} \subseteq {\mathcal P}(M)$ 
is {\it countably generated} if there is a countable ${\mathfrak X}_0\subseteq {\mathfrak X}$ 
that generates ${\mathfrak X}$. For example, $\Def(\MM)$ is countably generated 
since it is generated by the set of all those $X \in \Def(\MM)$ definable 
without parameters.

\bigskip


{\bf \S2. The Proof, I.} The purpose of this section is to prove Corollary~2.3, which, in the presence of Theorem~1, is the $(b) \Longrightarrow (a)$ half of Theorem~3. This is exactly what was proved in an earlier version of this paper, but here we 
 will prove the even stronger Theorem~2.1 in response to   a question raised by Roman Kossak. 
 
 \bigskip

Suppose that $\MM \prec \NN$. Recall that $\Lt(\NN / \MM)$ is the lattice of all $\KK$ such that 
$\MM \preccurlyeq \KK \preccurlyeq \NN$. We say that ${\mathbf D}$ is {\it dense} (for this extension) if 
${\mathbf D} \subseteq \Lt(\NN / \MM) \backslash \{\MM\}$ such that 
  whenever $\MM \prec \KK_0 \preccurlyeq \NN$, 
 then  there is $\KK_1 \in {\mathbf D}$ such that $\KK_1 \preccurlyeq \KK_0$. Let $\den(\NN / \MM)$ be the least cardinal $\delta$ such that $\delta = |{\mathbf D}|$ 
for some dense ${\mathbf D}$. 

\bigskip

{\sc Theorem 2.1}: {\em Suppose that $\MM \prec_{\sf end} \NN$ is countably generated and that $\delta = \den(\NN / \MM)$.   Then 
 every set that is $\Pi^0_1$-definable  in $(\MM, {\mathfrak X})$ is the union of  $\delta+ \aleph_0$ sets each of which is $\Sigma_1^0$-definable in $(\MM, {\mathfrak X})$.}

\bigskip

{\it Proof}.  Since $\NN$ is a countably generated extension of $\MM$, there is $\NN' \succ_{\sf end} \NN$ and $\NN'$ is a minimal elementary extension of $\MM$.  (The proof of this is a straightforward modification  of the proof of \cite[Theorem~2.1.12]{ksbook} concerning \mbox{superminimal} extensions.) Then 
$\cod(\NN' / \MM) = \cod(\NN / \MM)$ and $\den(\NN'/\MM) = \den(\NN / \MM)$. Thus, by replacing $\NN$ by $\NN'$, if necessary, we can assume that $\NN$ is a finitely generated extension of $\MM$. 
Let $c$ generate $\NN$ over $\MM$. 

Let ${\mathbf D}$ be a dense set such that $|{\mathbf D}| = \delta$. For every $\KK_0$ such that 
$\MM \prec \KK_0 \preccurlyeq \NN$, there is $\KK_1$ that is a finitely generated extension of $\MM$ such that $\MM \prec \KK_1 \preccurlyeq \KK_0$. Thus, we can safely assume that ${\mathbf D}$ consists only of finitely generated extensions of $\MM$. Let $\{c_\alpha : \alpha < \delta\} \subseteq N \backslash M$ be such that each $\KK \in {\mathbf D}$ is generated over $\MM$ 
by some $c_\alpha$. 

 For each $\alpha < \delta$, let $g_\alpha : M \into M$ be an $\MM$-definable function 
 such that $g_\alpha^\NN(c) = c_\alpha$. 

If $f,g : M \into M$ are $\MM$-definable functions and $A \in \Def(\MM)$, then we say that 
$f$ {\it refines} $g$ on $A$ if whenever $x,y \in A$ and $f(x) = f(y)$, then $g(x) = g(y)$. The  following claim 
is the purpose for introducing this definition.

\smallskip

{\sf Claim}: Suppose that $f : M \into M$ is an $\MM$-definable function. 
Then,  $f^\NN(c) > M$ iff there are $A \in \Def(\MM)$ and $\alpha < \delta $ such that $c \in A^\NN$ 
and $f$ refines $g_\alpha$ on $A$.

\smallskip 

We prove the Claim. Let $a = f^\NN(c)$.

First, suppose that $a > M$. Let $\alpha < \delta$ be such that $c_\alpha$ is in the elementary extension of $\MM$ generated by $a$.  Thus, there is an $\MM$-definable function  
$h : M \into M$ such that $h^\NN(a) = c_\alpha$. Let $A \in \Def(\MM)$ be the set of those $x \in M$ such that $hf(x) = g_\alpha(x)$. Then $f$ refines $g_\alpha$ on $A$, and 
since $h^\NN f^\NN(c) = c_\alpha = g_\alpha^\NN(c)$, then $c \in A^\NN$.

Next, suppose that 
 $A \in \Def(\MM)$ and $\alpha < \delta$ are such that  
$c \in A^\NN$ and $f$ refines $g_\alpha$ on $A$. Define $h : f[A] \into g_\alpha[A]$  so that if $x \in A$, then 
$hf(x) = g_\alpha(x)$. This defines $h$ uniquely, and $h$ is $\MM$-definable. Then, $h^\NN f^\NN(c) = g_\alpha^\NN(c)$  so that $h^\NN(a) = c_i$. Since $c_i > M$, then $f^\NN(c) = a > M$.

This completes the proof of the Claim.

\smallskip 

Now let $D \subseteq M$ be $\Pi_1^0$-definable in $(\MM, {\mathfrak X})$. We will prove that there 
are $\Sigma_1^0$-definable sets $D_{\alpha,j}$ ($\alpha < \delta$, $j < \omega$) whose union is~$D$. Since $D$ is $\Pi^0_1$-definable, we let $A \in {\mathfrak X}$ and  the ${\mathcal L}_\pa$-formula $\theta(x,y)$  be such that 
$$
D = \{d \in M : (\MM, {\mathfrak X}) \models \forall y \in A\, \theta(d,y)\}.
$$
To see this is possible, first observe that there is $B \in {\mathfrak X}$ such that 
$D = \{d \in M : (\MM, {\mathfrak X}) \models \forall z [\langle d,z \rangle \in B]\}$, and then let $A = M \backslash B$ and let $\theta(x,y)$ be $\forall z(y \neq \langle x,z \rangle)$.  
 Let $S$ be the branch of $2^{<M}$ that indicates $A$. Then $S \in {\mathfrak X} = \cod(\NN / \MM)$,  and $D$ is $\Pi_1^0$-definable from $S$ by the formula 
$$
\forall s \in S\forall k < \ell(s)\big[(s)_k = 0 \into \theta(x,(s)_k)\big].
$$ 
Since $S \in \cod(\NN / \MM)$, we assume that $c \in N \backslash M$ is such that 
$$S = \{ s \in M : \NN \models s \lhd c\}.$$

For each $d \in M$, define the function $f_d : 2^{<M} \into M$ so that{\footnote{This definition uses the minimalization operator $\mu$. Thus, if there is $k < \ell(s)$ such that $\MM \models  (s)_k = 0 \wedge \neg \theta(d,k)$, then $f_d(s)$ is the least such $k$; otherwise, $f_d(k) = \ell(s)$.}}
$$
f_d(s) = \mu k < \ell(s)[(s)_k = 0 \wedge \neg \theta(d,k)].
$$
Clearly, each $f_d$ is $\MM$-definable; in fact $\langle f_d : d \in M \rangle$ is an $\MM$-definable family of functions. It is easily seen that for each $d \in M$,
\begin{equation}\tag{2.1.1}
 d \in D  \Longleftrightarrow f_d^\NN(c) > M.
\end{equation}
From the Claim, we have that  $d \in D$ iff there are $A \in \Def(\MM)$ and $\alpha < \delta$ such that $c \in A^\NN$   and $f_d$ refines $g_\alpha$ on $A$. 
 
Let $\theta_0(x,y), \theta_1(x,y), \theta_2(x,y), \ldots$ be an enumeration  of all the $2$-ary ${\mathcal L}_\pa$-formulas. If $j < \omega$ and $b \in M$, let 
$A_j(b) \in \Def(\MM)$ be the set defined by $\theta_j(x,b)$. Then let 
$$
B_j = \{b \in M : c \in A_j^\NN(b)\},
$$
so that $B_j \in \cod(\NN / \MM) = {\mathfrak X}$. For $\alpha < \delta$ and $j < \omega$, let 
 $$
 D_{\alpha,j} = \{d \in M : {\mbox{ there is }} b \in B_j   \mbox{ such that} f_d {\mbox{ refines }} g_\alpha\mbox{ on } A_j(b)\}.
 $$
   Thus, each $D_{\alpha,j}$ is $\Sigma_1^0$-definable in $(\MM, {\mathfrak X})$. 
   There are at most $\delta + \aleph_0$ such $D_{\alpha,j}$. 
 
 To finish the proof, we will show that 
 $D = \bigcup\{D_{\alpha,j} : \alpha < \delta, \ j < \omega\}$.

 \smallskip

 $D \subseteq \bigcup\{D_{\alpha,j} : \alpha < \delta, \ j < \omega\}$: Let $d \in D$. Then $f_d^\NN(c) > M$ by (2.1.1). By the Claim, we let $\alpha < \delta$  and $A \in \Def(\MM)$ be such that $c \in A^\NN$ and $f_d$ refines  $g_\alpha$ on $A$. Let $j < \omega$ and $b \in B_j$ be such that $A = A_j(b)$. 
  Then,   $d \in D_{\alpha,j}$.  
  
  \smallskip
 
 $\bigcup\{D_{\alpha,j} : \alpha < \delta, \ j < \omega\} \subseteq D$: Suppose that $\alpha < \delta$, $j < \omega$ and $d \in D_{\alpha,j}$. Then there is $b \in B_j$ such that $f_d$ refines $g_i$ on $A_j(b)$. Since  $A_j(b) \in \Def(\MM)$ and $c \in A_j(b)^\NN$, 
 then $d \in D$, thereby completing the 
 proof of the theorem. \qed
 
 \bigskip
 
 {\sc Corollary 2.2}: {\em Suppose that $\MM \prec_{\sf end} \NN$ and $\Lt(\NN / \MM)$ is countable.    Then 
 every set that is $\Pi^0_1$-definable  in $(\MM, {\mathfrak X})$ is the union of countably many sets that are $\Sigma_1^0$-definable in $(\MM, {\mathfrak X})$.} \qed
 
 \bigskip
 
 The special case of the previous corollary when $|\Lt(\NN / \MM)| = 2$ should be mentioned.
 
 \bigskip
 
 {\sc Corollary 2.3}: {\em Suppose that $\MM \prec_{\sf end} \NN$ and $\MM$ is a minimal extension of $\NN$.      Then 
 every set that is $\Pi^0_1$-definable  in $(\MM, {\mathfrak X})$ is the union of countably many sets that are $\Sigma_1^0$-definable in $(\MM, {\mathfrak X})$.} \qed

\bigskip


\bigskip

{\bf \S3.\@ Outlines.} In this section we give an outline of the proof of the implication $(a) \Longrightarrow (c)$ of Theorem~1, and also of the implication that $(a)$ of Theorem~2 implies 
the weakening of $(c)$ in Theorem~2 obtained by replacing ``superminimal" with ``minimal''.  
We first consider  Theorem~1.

Suppose that $(\MM, {\mathfrak X})$ satisfies $(a)$ of Theorem~1; specifically, ${\mathfrak X}$ is countably generated, $\Def(\MM) \subseteq {\mathfrak X}$ and $(\MM, {\mathfrak X}) \models \wklsd$. We describe how $\NN$  as in $(c)$ 
was obtained in \cite{code}. A set $\Phi(x)$  of \mbox{$1$-ary} ${\mathcal L}_\pa(M)$-formulas is \mbox{{\it allowable}} if, for some $n < \omega$,  there are ${\mathcal L}_\pa$-formulas $\theta_0(x,y), \theta_1(x,y), $ $
\ldots, \theta_{n-1}(x,y)$ and $A_0,A_1, \ldots, A_{n-1} \in {\mathfrak X}$ such that 
$\Phi(x) = \bigcup_{j<n}\{\theta_j(x,b) : b \in A_j\}$ and whenever $\Psi(x) \subseteq \Phi(x)$ 
is $\MM$-finite, then $\bigwedge \Psi(x)$ defines an unbounded subset of $M$. 
Since $\varnothing \in {\mathfrak X}$, the set $\varnothing$ is allowable.

Our goal is to obtain an increasing sequence 
$$
\Phi_0(x) \subseteq \Phi_1(x) \subseteq \Phi_2(x) \subseteq \cdots
$$
 of allowable sets  that has the following two properties:

\begin{itemize}

\item[(A1)] (\cite[(4.1) and (5.1)]{code}) {\em Whenever $\theta(x,y)$ is an ${\mathcal L}_\pa$-formula, then there are $m < \omega$ and  $B \in {\mathfrak X}$ such that 
$$\{\theta(x,b) : b \in B\} \cup \{\neg\theta(x,b) : b \in M\backslash B\} \subseteq \Phi_m(x).$$}

\item[(A2)] (\cite[(4.2) and (5.2)]{code}){\em Whenever $B \in {\mathfrak X}$, then there are $m < \omega$ and an ${\mathcal L}_\pa$-formula $\theta(x,y)$ 
such that 
$$\{\theta(x,b) : b \in B\} \cup \{\neg\theta(x,b) : b \in M\backslash B\} \subseteq \Phi_m(x).$$}

\end{itemize}

If we have such a sequence, then $\bigcup_{m<\omega}\Phi_m(x)$ generates a complete type over $\MM$, thanks to (A1).
Let $\NN$ be an extension of $\MM$ generated by an element  realizing this complete type.
Then $\NN \succ_{\sf end} \MM$. It follows from (A1) that $\cod(\NN / \MM) \subseteq {\mathfrak X}$ and from (A2) that  $ {\mathfrak X}  \subseteq \cod(\NN / \MM) $.

The sequence is constructed by recursion, starting with $\Phi_0(x) = \varnothing$. 
Associated with  (A1)  is a lemma whose repeated application (once for each $\theta(x,y)$) assures that the constructed sequence of allowable sets satisfies (A1).  There is also such a lemma associated with (A2), which, when applied repeatedly (once for each $B$ in some countable set that generates ${\mathfrak X}$), assures that (A2) is satisfied. 
In each of these lemmas, we start with an allowable set $\Phi(x)$ which we then enlarge it to an allowable set $\Phi(x) \cup \{\theta(x,b) : b \in B\} \cup \{\neg\theta(x,b) : b \in M\backslash B\}$. 

The lemma for (A1) is the following.

\bigskip

{\sc Lemma 3.1}: {\em Suppose that $\Phi(x)$ is allowable  and that  $\theta(x,y)$ is 
an ${\mathcal L}_\pa$-formula. Then there is $B \in {\mathfrak X}$ such that 
$$
\Phi(x) \cup \{\theta(x,b) : b \in B\} \cup \{\neg\theta(x,b) : b \in M\backslash B\} 
$$
is allowable.}

\bigskip

The lemma for (A2) is the following.

\bigskip

{\sc Lemma 3.2:} {\em Suppose that $\Phi(x)$ is allowable  and   $B \in {\mathfrak X}$. 
Then there is an ${\mathcal L}_\pa$-formula $\theta(x,y)$ such that 
$$
\Phi(x) \cup \{\theta(x,b) : b \in B\} \cup \{\neg\theta(x,b) : b \in M\backslash B\} 
$$
is allowable.}

\bigskip

With  Lemmas~3.1 and~3.2, we are able to prove  $(a) \Longrightarrow (c)$ of Theorem~1 as in \cite{code}.

\bigskip

Next, we consider Theorem~2. In order to get $\NN \succ_{\sf end} \MM$ to be minimal, we require that the sequence of allowable sets satisfies the following additional condition.

\begin{itemize}

\item[(A3)] (\cite[(5.3)]{code}) {\em Whenever $f : M \into M$ is an $\MM$-definable function, then there are $m < \omega$ and $\varphi(x) \in \Phi_m(x)$ such that $f$ is either bounded or one-to-one on the set defined by $\varphi(x)$.}

\end{itemize}

 It is easily seen (or consult the proof of \cite[Theorem~5]{code}) that (A3) suffices to guarantee that $\NN$ is a minimal extension of $\MM$.

The next lemma is used to obtain (A3). By repeated applications of this lemma (once for each $\MM$-definable $f : M \into M$), we can get 
the sequence of allowable sets to satisfy (A3) as long as $\MM$ is countable. This will prove    the  weakening of $(a) \Longrightarrow (c)$ of Theorem~2, with ``minimal'' replacing ``superminimal''.

\bigskip

{\sc Lemma 3.3}: {\em Suppose that $\Phi(x)$ is allowable  and   $f : M \into M$ is $\MM$-definable. Then there is a formula $\varphi(x)$ such that $\Phi(x) \cup\{\varphi(x)\}$ is allowable and $f$ is either bounded or one-to-one on the set defined by $\varphi(x)$.} 

\bigskip

We remark that the proofs of  Lemmas~3.1, 3.2  and~3.3  do not require that  ${\mathfrak X}$ be countably generated. Furthermore, the proof of Lemma~3.3 does not require that $\MM$ be countable.

\bigskip


{\bf \S4.\@ The Proof, II.}  This section is devoted to completing the proof Theorem~3, the $(b) \Longrightarrow (a)$ half following from Corollary~2.3.  

Let $\MM \models \pa$ and ${\mathfrak X} \subseteq {\mathcal P}(M)$.
 Suppose that $(a)$ of Theorem~3 holds. By Theorem~1, we have that ${\mathfrak X}$ is countably generated, $\Def(\MM) \subseteq {\mathfrak X}$ and $(\MM, {\mathfrak X}) \models \wklsd$. In addition, every $\Pi_1^0$-definable set is the countable union $\Sigma_1^0$-definable  sets.  We will construct an increasing  sequence  $\Phi_0(x) \subseteq \Phi_1(x) \subseteq \Phi_2(x) \subseteq \cdots$ 
 of allowable sets satisfying the three conditions (A1), (A2) and (A3) from \S3.  The challenge will be to 
 satisfy (A3).
 
 Following are two corollaries of Lemma~3.1.

 \bigskip
 
 {\sc Corollary~4.1}: {\em Suppose that $\Phi(x)$ is allowable  and that  $\varphi(x,y)$ is an ${\mathcal L}_\pa$-formula. Then there are $B \in {\mathfrak X}$, an ${\mathcal L}_\pa$-formula $\theta(x,y)$  such that  
 $$
  \Phi(x) \cup  \{\theta(x,b) : b \in B\} 
   $$ 
   is allowable and 
for each $k \in M$ there is $b \in B$ such that 
$$
\MM \models \forall y < k\big[\big(\forall x(\theta(x,b) \into \varphi(x,y)\big) \vee 
                                             \big(\forall x (\theta(x,b) \into \neg\varphi(x,y)\big)\big].
$$}

\bigskip

{\it Proof}: Apply Lemma~3.1 to get $A \in {\mathfrak X}$ such that 
$$
\Phi(x) \cup \{\varphi(x,a) : a \in A\} \cup \{\neg\varphi(x,a) : a \in M\backslash A\} 
$$
is allowable. Let $B$ be the branch of $2^{<M}$ that indicates $A$. Let $\theta(x,y)$ be the formula
$$
\forall k < \ell(y)[(y)_k = 0 \leftrightarrow \varphi(x,k)],
$$
 and then let
 $$
 \Phi'(x) = \Phi(x) \cup \{\theta(x,b) : b \in B\}.
 $$
 It is clear that $\Phi'(x)$ is allowable. Moreover, for each $k \in M$, if we let $b \in  B$ be such that 
 $\ell(b) = k$, then $b$ is as required. \qed
  
 \bigskip

{\sc Corollary 4.2}: {\em Suppose that $\Phi(x)$ is an allowable set and $t(u,x)$ is   a Skolem ${\mathcal L}_\pa$-term. Then there are $B \in {\mathfrak X}$ and an ${\mathcal L}_\pa$-formula $\theta(x,y)$ such that $\Phi(x) \cup \{\theta(x,b) \in B\}$ is  allowable and  for each $k \in M$, there is $b \in B$  such that for $i,j \leq k$, either
\begin{equation}\tag{4.2.1}
\MM \models \forall x\big[ \theta(x,b) \into t(i,x) \leq  t(j,x) \big]
\end{equation}
or 
\begin{equation}\tag{4.2.2}
\MM \models \forall x\big[ \theta(x,b) \into t(i,x) >  t(j,x) \big].
\end{equation}}

\bigskip

{\it Proof}. Let $\varphi(x,y)$ be the formula $\exists i,j\big[y = \langle i,j \rangle \wedge t(i,x) \leq t(j,x)\big]$. Apply Corollary~4.1. \qed

 \bigskip
 
 The next lemma is a generalization of  Lemma~3.3.

\bigskip

{\sc Lemma 4.3}: {\em Suppose that $\Phi(x)$ is allowable, $t(u,x)$ is a Skolem 
${\mathcal L}_\pa$-term and $D \subseteq M$ is $\Sigma_1^0$-definable in $(\MM, {\mathfrak X})$. Suppose further that  whenever $i \in D$ and $a \in M$, then there is $\varphi(x) \in \Phi(x)$ such that 
$$\
\MM \models \forall x\big[\varphi(x) \into t(i,x) > a\big].
$$
 Then there is an allowable $\Phi'(x) \supseteq \Phi(x)$ such that 
whenever $i \in D$, there is $\varphi(x) \in \Phi'(x)$ such that 
$$
\MM \models \forall x,y\big[\varphi(x) \wedge \varphi(y) \wedge x < y \into t(i,x) < t(i,y)\big].
$$} 

\bigskip

{\it Proof}. Recall a definition from \cite{code}. An allowable set $\Psi(x)$ is  {\it  tree-based} if there is an ${\mathcal L}_\pa$-formula $\psi(x,y)$ and a branch $B \in {\mathfrak X}$ of the full binary tree $(M,\lhd)$ such that 
$$
\Psi(x) = \{\theta(x,b) : b \in B\}
$$ 
and 
$$
\MM \models \forall s,t\big[s \lhdeq t \into \forall x\big(\theta(x,t) \into \theta(x,s)\big)].
$$
Following \cite[p.\@ 576]{code}, we suppose, without loss of generality,  that $\Phi(x)$ is tree-based. Let $\Phi(x) = \{\theta(x,s) : s \in B\}$, where $B$ is a branch.  Let $T$ be the set of those $s \in 2^{<M}$ such that 
$\MM \models \forall w \exists x[x > w \wedge  \theta(x,s)]$. Then $T$ is $\MM$-definable without parameters. Also, $T$ is an $\MM$-tree and $B$ is a branch of $T$. 

For $i \in M$, let $f_i : M \into M$ be the function defined by $t(i,x)$. (It will be easier, notationally, to work with $f_i$ rather than $t(i,x)$.)

By Corollary~4.2, we can assume that whenever $k \in M$, then there is $s \in B$ such that whenever $i,j \leq k$, then either (4.2.1) or (4.2.2). Let $\gamma \subseteq M^2$ be such that 
whenever $i,j \in M$, then $\langle i,j \rangle \in \gamma$ iff 
\begin{equation}\tag{4.3.1}
(\MM, {\mathfrak X}) \models \exists s \in B \forall x\big[ \theta(x,s) \into f_i(x) \leq  f_j(x) \big].
\end{equation}
Equivalently, $\langle i,j \rangle \in \gamma$ iff 
\begin{equation}\tag{4.3.2}
(\MM, {\mathfrak X}) \models \forall s \in B \exists x\big[ \theta(x,s) \wedge f_i(x) \leq  f_j(x) \big].
\end{equation}
Hence, $\gamma$ is is $\Delta_1^0$-definable from $B$ and a set in $\Def(\MM)$, so that $\gamma \in {\mathfrak X}$. 

It follows from (4.3.1) that $\gamma$ is transitive, and from (4.3.2) that $\gamma$ is connected. Being both a transitive and connected binary relation on $M$, $\gamma$ is, by definition, a {\it linear quasi-order}  of $M$.  Notice that there is the possibility that $\langle i,j \rangle, \langle j,i \rangle \in \gamma$ for distinct $i,j \in M$, although  $\langle i,i \rangle \in \gamma$ for all $i \in M$.

Since $D$ is $\Sigma_1^0$-definable in $(\MM, {\mathfrak X})$, we let $A \in {\mathfrak X}$ and  the ${\mathcal L}_\pa$-formula 
$\psi(x,y)$ 
 be such that the formula $\exists y \in A\psi(x,y)$ defines $D$ in $(\MM, {\mathfrak X})$. 

Let $P = T \ \otimes \ 2^{<M}$ be the set of all pairs $\langle s,t \rangle$, where $s \in T, t \in 2^{<M}$ and 
$\ell(s) = \ell(t)$. We consider $P$ to be ordered componentwise by $\leq_P$; that is, 
if $p = \langle s_0,t_0 \rangle$ and $q =  \langle s_1,t_1 \rangle$ are in $P$, then 
$p \leq_P q$ iff $s_0 \lhdeq s_1$ and $t_0 \lhdeq t_1$. 
Then $(P, \leq_P)$ is (isomorphic to) an $\MM$-tree that is definable in $\MM$ without parameters. 
Let 
$$
C = \{\langle s,t \rangle \in P : s \in B {\mbox{ and }} \MM \models \forall i < \ell(t)[(t)_i = 0 \leftrightarrow i \in A]\}.
$$ 
 Then, $C \in {\mathfrak X}$ and $C$ is a branch of $P$. 
 
 \smallskip

{\it Convention}: Even though $P \subseteq M^2$, by  the usual coding of pairs, 
we also assume that  $P 
\subseteq M$. We then have that $<$ linearly orders $P$. 
We assume the coding is done so that whenever $q <_P p \in P$, then $q<p$.

\smallskip

 We will define $P_0 \subseteq P$ and a function $g : P_0 \into M$ by recursion in~$\MM$. 
 Suppose that $p \in P$ and that, for all $q < p$, we have already decided whether or not $q$ is in $ P_0$ and, if it is, what $g(q)$ is. Then,  $p \in P_0$ iff  there is $u \in M$ such that both:
 
 \smallskip
 
  \begin{itemize}

\item[(4.3.3)] \quad $\MM \models \theta(u,s) \wedge \forall q \big[\big(q \in P_0 \wedge q < p\big)  \into u > g(q)\big]$;

\item[(4.3.4)] \quad
$\MM \models \forall i,q\Big[\Big(  q \in P_0 \wedge q < p \wedge  i < \ell(s) \wedge\\ {\mbox{   }}  \quad \quad \quad \exists k < \ell(t)\big((t)_k = 0 \wedge \psi(i,(t)_k)\big)\Big) \into f_i(u) > f_i(g(q))\Big]$.

\end{itemize}
If $p \in P_0$,  let $g(p)$ be the least  $u$ such that $(4.3.3)$ and $(4.3.4)$ hold. 

Both $P_0$ and $g$ are $\MM$-definable without parameters.

If $q <_P p  \in P_0$, then $q \in P_0$. For if $u \in M$ is a witness 
demonstrating  that $p \in P_0$ (that is, $u$ satisfies both (4.3.3) and (4.3.4)), then it already had demonstrated that $q \in P_0$.
Thus, we have that $P_0$ is a subtree of~$P$. 

We claim that $C \subseteq P_0$. Consider $p = \langle s,t \rangle \in C$, intending to show that 
$p \in P_0$.  Thus, we want $u \in M$ that satisfies (4.3.3) and (4.3.4).   Since $s \in B$, then $\theta(x,s) \in \Phi(x)$ and, therefore, $\theta(x,s)$ defines an unbounded subset of $M$. Hence,  there is an unbounded set of $u$'s satisfying (4.3.3). 

For (4.3.4), the $f_i$'s that need to be considered are those for which 
$$
\MM \models i < \ell(s) \wedge \exists k \big[k < \ell(t) \wedge \psi(i,(t)_k)\big].
$$
 Let $I$ be the set of such $i$. Clearly, $I \subseteq D$. Let $j \in I$ be  a  $\gamma$-minimal element of $I$; that is,   $\langle j,i \rangle \in \gamma$ for all $i \in I$. (Since $I, \gamma \in {\mathfrak X}$ and $I$ is bounded, then  such a $j$ does exist.)
 
Let 
$$
a = \max\{f_i(g(q)) : p > q \in P_0, i < \ell(s)\}.
$$ 
Let $s' \in B$ be large enough so that $s' \rhd s$, 
\begin{equation}\tag{4.3.5}
\MM \models \forall x\big[\theta(x,s') \into f_{j}(x) > a \big]
\end{equation}
and
\begin{equation}\tag{4.3.6}
\MM \models \forall i,j \in I \forall x\big[\theta(x,s') \into \big(f_i(x) \leq f_j(x) \leftrightarrow \langle i,j \rangle \in \gamma \big) \big].
\end{equation}
For (4.3.5), $s'$ exists by the hypothesis of the lemma since $I \subseteq D$. 
For (4.3.6),  invoke Corollary~4.1.

Then, for any $i \in I$, we have that 
$$
f_i(u) \geq f_{j}(u) > a \geq  f_i(g(q))
$$
whenever  $u$ and $q$ are such that $\MM \models \theta(u,s')$ and $p > q \in P_0$. Therefore, $p \in P_0$.

Let $\theta'(x,p)$ be the ${\mathcal L}_\pa$-formula 
$$
 \exists q \big(p \leq_P q \in P_0  \wedge g(q) = x\big)
$$ 
and  let $$
\Phi'(x) = \Phi(x) \cup \{\theta'(x,p) : p \in C\}.
$$  It is seen that $\Phi'(x)$ is allowable. 
We claim that $\Phi'(x)$ has the property required by the lemma.

To see that it does, consider $i \in D$. Let $p_0 = \langle s_0,t_0 \rangle \in C$ be large such  that  
$\ell(s_0) = i+1$. 
Let $\varphi(x) = \theta'(x,p_0) \in \Phi'(x)$, and let $x,y \in M$ be such that 
$\MM \models \varphi(x) \wedge \varphi(y) \wedge x < y$. Then there are $q,r \geq_P p_0$ 
such that $q,r \in P_0$ and $g(q) = x$ and $g(r) = y$. Then (4.3.3) implies that $q < r$. 
Let $r = \langle s, t \rangle$. Since $t = t_0 \harp (i+1)$, then $i$ satisfies the left side of the implication in (4.3.4), and therefore, $f_i(y) > f_i(z)$. This completes the proof of the lemma. \qed

\bigskip

We remark that the proofs of Corollary~4.1 and~4.2 and Lemma~4.3 do not require that ${\mathfrak X}$ be countably generated nor that every $\Pi^0_1$-definable set is the union of countably many 
$\Sigma_1^0$-definable sets.

We now construct the sequence $\Phi_0(x) \subseteq \Phi_1(x) \subseteq \Phi_2(x) \subseteq \cdots$. By repeated applications of Lemma~3.1 (once for each  $\varphi(x,y)$), we will get (A1) to hold. By repeated applications of Lemma~3.2 (once for each $B$ in some countable set of generators for ${\mathfrak X}$), we will get (A2) to hold. We now describe how we get (A3) to hold.  

Consider a Skolem ${\mathcal L}_\pa$-term $t(u,x)$. Let $\varphi(x,y)$ be the formula 
$\exists u,v\big[y = \langle u,v \rangle \wedge t(u,x) \leq v\big]$. At some point in the construction, we will have applied Lemma~3.1 to this formula, obtaining $\Phi_m(x)$. Therefore, there are 
$B \in {\mathfrak X}$ and an ${\mathcal L}_\pa$-formula $\theta(x,y)$ such that 
$\{\theta(x,b) : b \in B\} \subseteq \Phi_m(x)$ and for 
every $i,a \in M$, there is $b \in B$ such either 
\begin{equation}\tag{4.1}
\MM \models \forall x\big[\theta(x,b) \into t(x,i) \leq a\big]
\end{equation}
or
\begin{equation}\tag{4.2}
\MM \models \forall x\big[\theta(x,b) \into t(x,i) > a\big].
\end{equation}
Let $D$ be the set of all those $i \in M$ for which (4.1) holds for every $a \in M$. Thus, the set $D$ is $\Pi_1^0$ in $B$ and some set in $\Def(\MM)$. 
By the condition on $(\MM, {\mathfrak X})$, there are countably many $\Sigma_1^0$-definable  sets 
$D_0,D_1,D_2, \ldots$ whose union is $D$. At any future stage of the construction, where we have (say) $\Phi_n(x)$, we will have, for any $j < \omega$, that the hypothesis of Corollary~4.1 holds (with $\Phi(x) = \Phi_n(x)$ and $D = D_j$). Applying Lemma~4.3, we get $\Phi_{n+1}(x) = \Phi'(x)$ so that whenever $i \in D_j$, then  there is $\varphi(x) \in \Phi_{n+1}(x)$ such that 
$$
\MM \models \forall x,y\big[\varphi(x) \wedge \varphi(y) \wedge x < y \into t(i,x) < t(i,y)\big].
$$
Thus by repeated applications of Lemma~4.3 (once for each  $D_j$), we will get that for every $i \in M$, there are $n < \omega$ and $\varphi(x) \in 
\Phi_n(x)$ such that either
$$
\MM \models \exists w\forall x\big[\varphi(x) \into t(i,x) \leq w\big]
$$
or
$$
\MM \models \forall x,y\big[\varphi(x) \wedge \varphi(y) \wedge x < y \into t(i,x) < t(i,y)\big].
$$ 
Since there are only countably many such $t(u,x)$ we can do the construction by dovetailing the steps necessary for each  possible $t(u,x)$. 

Thus, we get the sequence $\Phi_0(x) \subseteq \Phi_1(x) \subseteq \Phi_2(x) \subseteq \cdots$ to satisfy (A3), thereby completing the proof of the $(a) \Longrightarrow (b)$ half of Theorem~3. 

With Theorem~1 and Corollary~2.3, the proof of Theorem~3 is complete.  

\bigskip


{\bf \S5.\@ Appendix.} It may not be apparent that the added condition in (a) of Theorem~3 is actually needed. In other words, perhaps whenever  $\MM, \NN$ are such that $\NN \succ_{\sf end} \MM$ is a finitely generated extension, then $\MM$ has a minimal extension $\NN_0 \succ_{\sf end} \MM$ such that $\cod(\NN_0 / \MM) = \cod(\NN / \MM)$.
The next theorem shows that that is not so.

\bigskip

{\sc Theorem 5.1}: {\em For every $\MM_0 \models \pa$, there is $(\MM, {\mathfrak X}) \models \wkl $ 
such that $\MM \equiv \MM_0$, ${\mathfrak X}$ is countably generated, ${\mathfrak X} \supseteq \Def(\MM)$  and there is a $\Pi_1^0$-definable set that is not the union of countably many $\Sigma_1^0$-definable sets.}

\bigskip

{\it Proof}. Suppose that  $\MM_0 \models \pa$. Viewing $(\MM_0, \Def(\MM_0))$  
as a \mbox{2-sorted} first-order structure, we let 
$(\MM, {\mathfrak X}_0) \equiv (\MM_0, \Def(\MM_0))$ be  \mbox{$\aleph_1$-saturated}.   By $\aleph_1$-saturation (although only recursive saturation is needed), there is $X \in {\mathfrak X}_0$ 
 such that $(\MM, {\mathfrak X}_0)  \models (X)_n = \varnothing^{(n)}$ for every $n < \omega$.
(That is, $(X)_n$ is the complete $\Sigma_n^0$-definable subset of $\MM$.) 
Obtain by recursion, using the Low Basis Theorem,  a sequence $X_0,X_1,X_2, \ldots$ such that $X_0 = X$ and, for each $n < \omega$,  
$X_{n+1} \in {\mathfrak X}_0$ is low relative to $X_n$ (i.e., that $(\MM, {\mathfrak X}_0) \models (X_{n+1})' \equiv_T  (X_n)'$) and is such that for any unbounded 
$\MM_0$-tree $S$ that is $\Delta_1^0$-definable from $X_0 \oplus X_1 \oplus \cdots \oplus X_n$, there is $k \in M$ is such that $(X)_k$ is a branch of $S$. 
Let ${\mathfrak X}$ be the set of those subsets of $M$ that are $\Delta_1^0$-definable in  
$\Def((\MM,  X_0,X_1,X_2, \ldots))$. 
We easily see that $(\MM, {\mathfrak X})$ is as required.

Obviously, $\MM \equiv \MM_0$. The countable set $\{X_0,X_1,X_2, \ldots\}$ 
generates~${\mathfrak X}$.
Since each $X_{n+1}$ is low relative to $X_n$, then each $X_n$ is low relative to $X$, so that $X' \not\in {\mathfrak X}$. Thus, $M \backslash X'$ is a $\Pi_1^0$-definable set that is not $\Sigma_1^0$-definable, nor is it the union of finitely many $\Sigma_1^0$-definable sets. By $\aleph_1$-saturation, it is not the union of countably many $\Sigma_1^0$-definable sets.
\qed

\bigskip

{\sc Remark}: In the proof of  the previous theorem, if $(\MM, {\mathfrak X}_0)$ happened to be 
$\kappa^+$-saturated, where $\kappa \geq \aleph_0$,  then  $M \backslash X'$ would be a $\Pi^0_1$-definable set that is not the union 
of $\kappa$ $\Sigma^0_1$-definable sets.

\bibliographystyle{plain}

\end{document}